\documentclass[journal, 10pt]{IEEEtran}
\usepackage{amsfonts,amsmath,amssymb,amsthm}
\usepackage{algorithmic}
\usepackage{algorithm}
\usepackage{times}
\usepackage{bm}
\usepackage{graphicx}
\usepackage{setspace}
\usepackage{cases}
\usepackage{comment}
\graphicspath{{./fig_eps/}}
\usepackage{here}
\usepackage{ascmac} 
\usepackage{array}
\usepackage{float}
\usepackage{cite}
\usepackage{url}
\usepackage{subfigure}
\usepackage{ulem}

\theoremstyle{remark}
\newtheorem{mydef}{Definition}
\newtheorem{mylem}{Lemma}
\newtheorem{mythm}{Theorem}
\newtheorem{myas}{Assumption}
\newtheorem{myrem}{Remark}

\newcommand{\rfig}[1]{Fig.~\ref{#1}} 
\newcommand{\req}[1]{\eqref{#1}} 
 
\newcommand{\rtab}[1]{Table~\ref{#1}}
\newcommand{\rlem}[1]{Lemma~\ref{#1}}
\newcommand{\ras}[1]{Assumption~\ref{#1}}
\newcommand{\rrem}[1]{Remark~\ref{#1}}
\newcommand{\rdef}[1]{Definition~\ref{#1}}

\begin{document}
\title{Self-triggered Model Predictive Control for Nonlinear Input-Affine Dynamical Systems \\via Adaptive Control Samples Selection}
\author{\large Kazumune~Hashimoto, Shuichi~Adachi,~\IEEEmembership{Member,~IEEE,} and~Dimos~V.~Dimarogonas,~\IEEEmembership{Member,~IEEE}
\thanks{Manuscript submitted to IEEE Transaction on Automatic Control.}
\thanks{Kazumune Hashimoto is with the Department of Applied Physics and Physico-Informatics, Keio University, Yokohama, Japan (e-mail: kazumune.hashimoto@z5.keio.jp).}
\thanks{Shuichi Adachi is with the Department of Applied Physics and Physico-Informatics, Keio University, Yokohama, Japan (e-mail: adachi@appi.keio.ac.jp).}
\thanks{Dimos V. Dimarogonas is with the ACCESS Linnaeus Center,
School of Electrical Engineering, KTH Royal Institute of Technology, 10044 Stockholm, Sweden (e-mail : dimos@ee.kth.se). His work was supported by the Swedish Research council (VR) and the Knut and Alice Wallenberg Foundation.}
}
\date{}
\maketitle
\begin{abstract}
In this paper, we propose a self-triggered formulation of Model Predictive Control for continuous-time nonlinear input-affine networked control systems. Our control method specifies not only when to execute control tasks but also provides a way to discretize the optimal control trajectory into several control samples, so that the reduction of communication load will be obtained. 
Stability analysis under the sample-and-hold implementation is also given, which guarantees that the state converges to a terminal region where the system can be stabilized by a local state feedback controller. 
Some simulation examples validate our proposed framework.
\end{abstract}
\section{Introduction}
\IEEEPARstart{E}{vent-triggered} control is one of the sampled-data control schemes that has been receiving increased attention in recent years \cite{ref1, ref2, ref3, ref4, ref5, ref6, ref7, ref8, ref9, ref10, ref11, ref12, ref13, ref14, ref15, ref16, ref17, ref18, ref19, ref20, ref21, ref22, ref23, ref23_2, ref24}. In contrast to time-triggered control where the control execution is periodic, event-triggered control requires the executions only when desired control specifications cannot be guaranteed. This may have several advantages over time-triggered control for  networked control systems, since this leads to the reduction of over-usage of communication resources and energy consumption when limited battery powered devices exist. Two main event-triggered control approaches have been proposed, namely event-based control \cite{ref1, ref2, ref3, ref4, ref5, ref6, ref7, ref8}, and self-triggered control \cite{ref9, ref10, ref11, ref12}. 
The main difference of these two approaches is that, for the event-based case the control input is executed based on the continuous state measurement of the plant, while for the self-triggered case the control execution is pre-determined based on the prediction from the plant model.

The event-triggered control framework has been analyzed for many different types of systems with different performance guarantees. 
For example, in \cite{ref8}, \cite{ref9} the authors propose the event-triggered strategy based on ${\cal L}_2$ and ${\cal L}_{\infty}$ gain stability performance for linear systems. Another research line formulates the triggering rules based on Input-to-State Stability (ISS) for linear systems, e.g., \cite{ref10}, which is followed by the extension to the nonlinear case \cite{ref18}. 

In this paper, we are interested in applying the event-triggered scheme to Model Predictive Control (MPC). 
This has been motivated due to the fact that MPC not only takes into account several constraints such as actuator limitations explicitly by solving the Optimal Control Problem (OCP) on-line, but also stability can be analyzed even for the nonlinear case. 
The application of event-triggered strategies to MPC has been receiving attention in recent years and some results include \cite{ref13, ref14, ref15, ref16, ref17, ref18, ref19, ref20, ref21, ref22, ref23, ref23_2}, where many results have been proposed for discrete time systems.  
For example, in \cite{ref13}, the authors propose event-based strategy for discrete time linear system under additive bounded disturbances, where the OCP is solved when the difference of actual state and predictive state exceeds a certain threshold, and they also analyze stability depending on the size of the disturbances. 
The reader can also refer to \cite{ref15}, and to more recent results in \cite{ref16}, \cite{ref17} for linear systems where infinite horizon quadratic cost is evaluated. The concept of event-triggered MPC has also been applied to cooperative control of multi-agent systems, see e.g., \cite{ref18, ref19}. For the continuous case, the reader can refer to \cite{ref20}, \cite{ref21}, \cite{ref22}. For example, in \cite{ref20}, the authors derive self-triggered MPC based on the optimal cost function as a Lyapunov candidate for controlling a non-holonomic vehicle, where the controlled plant solves an OCP only when it is needed. 

The contribution of this paper is to propose a new self-triggered MPC framework for nonlinear input-affine continuous-time networked control systems, where the plant with actuator and sensor systems are connected to the controller through wired or wireless channels. The derivation of our self-triggered strategy follows the previous ideas that the OCP is solved by checking if the optimal cost function, regarded as a Lyapunov candidate, is decreasing. 
In networked control systems, however, one of the main constraints is the bandwidth limitation \cite{ref25}, meaning that the controller can send only a limited number of control samples. In the up-to-date results of event-triggered MPC for continuous systems presented in \cite{ref20}, \cite{ref21}, \cite{ref22}, the event-triggered strategies for the continuous system are based on the assumption that the current and future (continuous) optimal control trajectory can be applied until the next OCP is solved. Therefore, this framework cannot be applied to our case, since continuous information cannot be transmitted to the plant needing an infinite transmission bandwidth. Even though the continuous control trajectory may be approximated by using a large number of control samples, it may lead to the over-usage of communication resources since it requires transmissions of many control samples.

Therefore, in our proposed control method, the controller not only solves an OCP but also \textit{discretizes} the obtained optimal control input trajectory into several control input samples, so that these can be transmitted as a packet to the plant. Then, the plant applies the control samples as in a sample-and-hold implementation. 
The discretizing method is to some extent relevant to ``Roll-out event-triggered control'', which is introduced in \cite{ref17}, where the authors propose a way to pick up the transmission time step for linear discrete time systems, and then show that the proposed control policy provides better performance than the conventional periodic optimal control in terms of the reduced value function. 
In contrast to \cite{ref17}, we will propose a way to adaptively select sampling time intervals to reduce the communication load. While this may lead to additional optimization problems, 
we will provide an efficient way of choosing the sampling intervals. 
Moreover, while the results presented in \cite{ref17} considers linear systems, 
we deal with nonlinear systems.

One of the main difficulties regarding MPC under sample-and-hold implementation is to guarantee stability, since sample-and-hold controllers lead to an error between the predicted optimal state and the actual state of the system, even when the system has no disturbances.
Regarding this stability problem, some results were provided in \cite{ref26},\cite{ref27}. The key idea is to use Lyapunov-based MPC, where a Lyapunov based controller is assumed to exist to show that the state converges to a certain invariant set under sample-and-hold implementation. 
However, it was not concluded whether the state converges to the terminal region where an assumed local controller exists stabilizing the system to the origin. 
In the MPC framework, it is desirable to achieve the convergence to the terminal region, since then the local state feedback controller can be applied to stabilize the system without needing to solve the OCP. The strategy of switching MPC to the local controller is referred to as `Dual-mode MPC'. Motivated by this, in this paper we also show that the state reaches the terminal region in finite time. Instead of using Lyapunov based MPC, an additional control input constraint and a restricted terminal constraint are used. 

As illustrative examples, we simulate both linear and nonlinear systems. For the linear system case, the problem of stabilizing an un-stable system under control input constraints is considered and we compare the control performance with periodic MPC under sample-and-hold implementation with the same average transmission interval. For the nonlinear system case, we consider position control of a non-holonomic vehicle in two dimensions. For both cases, we will show that the system is stabilized under our proposed self-triggered MPC.

The remainder of this paper is organized as follows. In Section II, the problem formulation is set up for the networked control system. In Section III, the self-triggered rule is given. In Section IV, we propose an efficient way to choose optimal sampling intervals. In Section V, 
the stability analysis is given. In Section VI, some simulation results are given. Finally, we summarize the results of this paper in Section VII.
\section{Problem Formulation}
\subsection{Notations}
We make use of the following notations.
Let $\mathbb{R}$,  $\mathbb{R}_{\geq 0}$,  $\mathbb{N}_{\geq 0}$, $\mathbb{N}_{\geq 1}$ be the real, non-negative real, non-negative integers and positive integers, respectively.
The operator $||\cdot ||$ denotes Euclidean norm of a vector.
The function $\phi (x,u): \mathbb{R}^n \times \mathbb{R}^m \rightarrow \mathbb{R}^n$ is locally Lipschitz continuous with Lipschitz constant $L_{\phi}$ in $x\in \Omega \subset \mathbb{R}^n $, if $||\phi (x_1, u) - \phi (x_2, u) || \leq L_{\phi} ||x_1 - x_2||$ for all $x_1, x_2 \in \Omega$. 
The difference between two sets $\Omega_1$ and $\Omega_2$ is denoted by $\Omega_1 \backslash \Omega_2 =\{ x\ |\ x\in \Omega_1, x\notin \Omega_2 \}$. A continuous function $\alpha$ : $[0, a) \rightarrow [0, \infty)$ is said to be ${\cal K}_{\infty}$ function if $\alpha$ is strictly increasing with $\alpha(0) =0$, and $\alpha(r) \rightarrow \infty$ as $r\rightarrow \infty$.
Given a compact set $\Omega \subseteq \mathbb{R}^n $, we denote $\partial \Omega$ as the boundary of $\Omega$.
\subsection{System Definition}
\begin{figure}[t]
  \begin{center}
   \includegraphics[width=9.0cm]{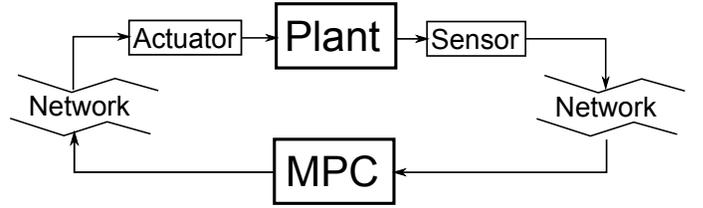}
   \caption{Networked Control System}
   \label{network}
  \end{center}
 \end{figure}

Consider the networked control system in \rfig{network}, where the plant with sensor and actuator systems are connected to the Model Predictive Controller (MPC) through the network channels. The dynamics of the plant are given by the following continuous-time nonlinear input affine system:
\begin{equation}\label{sys1}
\dot{{x}}(t) = \phi ({x} , u) = f({x}) + g({x}) u
\end{equation}
where ${x}\in \mathbb{R}^n$ is the state of the plant, and $u\in \mathbb{R}^m$ is the control input. We assume that the constraint for the control input is given by $||u|| \leq u_{\rm max}$. 
Our control objective is to asymptotically stabilize the system \req{sys1} to the origin, i.e., $x(t) \rightarrow 0$ as $t \rightarrow \infty$. To achieve this goal, we assume that the nonlinear system given by \req{sys1} satisfies the following conditions: 
\begin{myas}\label{as1}
The nonlinear function $\phi (x,u): \mathbb{R}^n \times \mathbb{R}^m  \rightarrow \mathbb{R}^n$ satisfies $\phi (0,0) = 0$, and is Lipschitz continuous in $x \in \mathbb{R}^n$ with Lipschitz constant $L_{\phi}$.
Furthermore, there exists a positive constant $L_G > 0$, such that $|| g(x) || \leq L_G$. 
\end{myas}

\begin{figure}[t]
  \begin{center}
   \includegraphics[width=7cm]{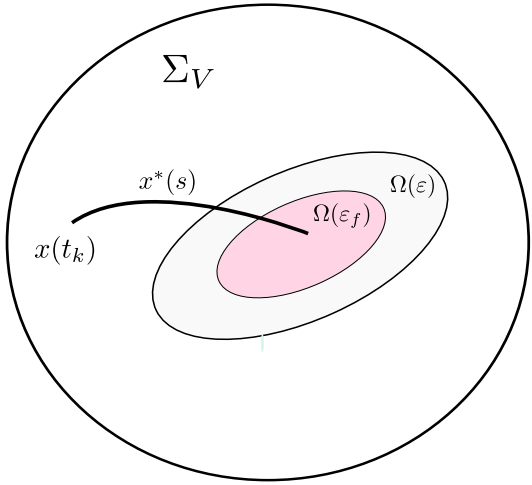}
   \caption{The illustration of three regions $\Sigma_V$, $\Omega (\varepsilon )$, $\Omega (\varepsilon_f )$, and an example of optimal state trajectory $x^* (s)$.}
   \label{ilterregion}
  \end{center}
 \end{figure}
\subsection{Optimal Control Problem}
In this subsection the Optimal Control Problem (OCP) is defined. 
Let the sequence $\{ t_k \}_{k\in \mathbb{N}_{\geq 0}}$ denote the sampling instants when the OCP is solved. At $t_k$, the controller solves the OCP involving the predictive states denoted as ${x}(s)$, and the control input $u(s)$ for $s\in [t_k, t_k +T_p]$ based on the current state measurement $x(t_k)$, where $T_p$ is the prediction horizon. In this paper, we consider the following cost function to be minimized (not necessarily quadratic costs)
\begin{equation*}
\begin{array}{lll}
J(x(t_k), u(\cdot )) ={\displaystyle \int}^{t_k+T_p}_{t_k}F({x}(s), u(s )) {\rm d}s+V_f({x}(t_k+T_p))
\end{array}
\end{equation*}
where $F(x, u)$ and $V_f(x)$ is a stage and terminal cost, and several assumed conditions are described later in this section.

Using the cost function above, the OCP to be solved is defined as follows:

\noindent 
\textit{Problem} 1 {(OCP) }: {At any update time $t_k$, $k \in \mathbb{N}_{\geq 0}$, given $x(t_k)$ and $T_p$, find the optimal control input and corresponding state trajectory $u^* (s )$, ${x}^* (s )$ for all $s \in [t_k , t_k +T_p ]$ that minimizes $J(x(t_k), u(\cdot ))$, subject to}
\begin{numcases}
   { }
      \dot{{x}}(s)=\phi ({x}(s) , u(s) ), \ s \in [t_k, t_k+T_p] \label{constraint1} \\ 
    {u}(s ) \in {\cal U} \label{constraint2} \\
    {x} (t_k + T_p )\in \Omega (\varepsilon_f) , \label{constraint3}
\end{numcases}
{ where ${\cal U}$ is the control input constraint set given by}
\begin{equation}\label{controlinputconstraint}
{\cal U}=\{ u(s) \in \mathbb{R}^m : ||u(s) ||\leq u_{\rm max}, ||\dot{u}(s)|| \leq K_u \}, 
\end{equation}
{and for given $\varepsilon _f > 0$, $\Omega (\varepsilon_f)$ is the terminal constraint set given by }
\begin{equation}
\Omega (\varepsilon_f) = \{ x \in \mathbb{R}^{n} : V_f (x) \leq \varepsilon _f \}.
\end{equation}

Regarding the control input constraint set in \req{controlinputconstraint}, we additionally consider a constant $K_u$ satisfying $||\dot{u}(s) || \leq K_u$. Although this puts a limit on the slope of the optimal control input and is sometimes used when the actuator has a physical limitation with the rate of its position change \cite{ref19}, \cite{ref28}, 
in this paper we will make use of this constraint to guarantee the stability analyzed in subsequent sections. 
For the (terminal) set $\Omega (\cdot )$, we further assume the following.
\begin{myas}\label{as2}
There exists a positive constant $\varepsilon>\varepsilon_f$ and a local stabilizing controller $\kappa(x) \in {\cal U}$, satisfying 
\begin{equation}\label{terconst}
\cfrac{\partial V_f}{{\partial}x} ( f(x) + g(x) \kappa(x) ) \leq - F(x, \kappa(x) ),  
\end{equation}
for all $x \in \Omega (\varepsilon )$.
\end{myas}
Note that since $\varepsilon > \varepsilon_f$, the terminal set $\Omega (\varepsilon_f )$ used in Problem 1 is smaller than $\Omega (\varepsilon )$ in \ras{as2}, i.e., $\Omega (\varepsilon_f )  \subset \Omega (\varepsilon )$. 

Denote $J^* (x(t_k))$ as the optimal cost obtained by Problem~1
\begin{equation*}
J^* (x(t_k)) = \underset{u(\cdot )}{\rm min}\ J(x(t_k), u(\cdot )).
\end{equation*}
Moreover, we consider a following set as a stability region characterized by $J^* (x(t_k))$ : 
\begin{mydef} \label{def2}
$\Sigma_V$ is the set given by $\Sigma_V = \{ x\in \mathbb{R}^n : J^* (x) \leq J_0 \}$, where $J_0$ is defined such that $\Omega (\varepsilon) \subseteq \Sigma_V$.
\end{mydef}
The illustration of the three regions considered in this paper $\Sigma_V$, $\Omega (\varepsilon)$, $\Omega (\varepsilon_f)$, and an example of optimal trajectory $x^* (s)$ that is constrained to be in $\Omega (\varepsilon_f)$ by the prediction horizon $T_p$, are all shown in \rfig{ilterregion}. 

We will show in this paper that if the state initially starts from inside the set $x \in \Sigma_V \backslash \Omega (\varepsilon) $, then the state trajectory enters $\Omega (\varepsilon)$ in finite time. Since the local control law $\kappa(x)$ is given from \ras{as2}, the system \req{sys1} can be stabilized by using $\kappa(x)$ once the state reaches $\Omega (\varepsilon)$ without needing to solve the OCP. For this reason, we consider that the control law switches from the solution to Problem 1 to the utilization of $\kappa(x)$ once the state enters $\Omega (\varepsilon)$.
This control scheme is in general referred to as `Dual-mode MPC', and is adopted in many works in the literature, see e.g., \cite{ref21}, \cite{ref33}. 
\begin{figure}[t]
  \begin{center}
   \includegraphics[width=7.0cm]{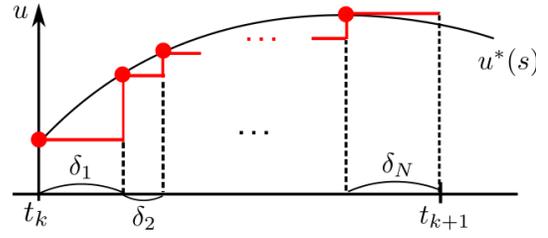}
   \caption{Optimal control input obtained at $t_k$: the controller picks up $N$ control input samples (Red circles) from the obtained optimal control trajectory (Black line), and these samples are transmitted to the plant and applies them as sample-and-hold fashion (Red line). }
   \label{controlinput}
  \end{center}
 \end{figure}
Using the sets defined above, the following conditions are further assumed to be satisfied for the stage and terminal cost $F$, $V_f$: 
\begin{myas}\label{as3}
$F(x,u) \geq \alpha_1 (||x||)$, and $V_f (x)\leq \alpha_2 (||x||) $, where $\alpha_1$ and $\alpha_2$ are ${\cal K}_{\infty}$ functions. Moreover, $F(x,u)$ and $V_f(x)$ are Lipschitz continuous in $x\in \Sigma_V$ with corresponding Lipschitz constants $0<L_F<\infty$, $0<L_{V_f}<\infty$.   
\end{myas}
\begin{myrem}
Assumptions 2 and 3 are fairly standard assumptions to guarantee stability for nonlinear systems under MPC. Several methods to numerically obtain $\Omega (\varepsilon)$ and $\kappa(x)$ satisfying \req{terconst} were proposed in e.g., \cite{ref32}, \cite{ref33}. 
Furthermore, several ways to compute Lipschitz parameters $L_F$, $L_{V_f}$ for the case of quadratic stage and terminal cost have been proposed in \cite{ref20}. 
\end{myrem}
In the following, let the optimal control input and the state trajectories obtained by Problem~1 be given by 
\begin{equation}\label{optux}
u^{*} (s ), \ \  {x}^{*}(s),\ \ s \in[t_k ,t_k + T_p]
\end{equation}
where $x^*(t_k) = x(t_k)$. 

Note that in earlier results of periodic or event-triggered MPC for continuous systems, e.g., \cite{ref20}, \cite{ref21},  \cite{ref22}, \cite{ref33}, \cite{ref36}, \cite{ref41}, the current and future continuous optimal control trajectory $u^*(s)$ is considered to be applied to the plant for $s\in[t_k, t_{k+1}]$. 
However, this situation cannot be applied to the networked control system in \rfig{network} considered in this paper, since sending continuous information requires an infinite transmission bandwidth. 
Therefore, we consider that only $N$ ($N\in \mathbb{N}_{\geq 1}$) control input samples, i.e., 
\begin{equation}\label{sampledcontroller}
\{ u^{*}(t_k), u^* (t_k +\delta_1), \cdots, u^* (t_k + \sum^N _{i=1} \delta_i ) \},
\end{equation}
should be determined to be picked up by the controller and then transmitted to the plant. The plant then applies the obtained control inputs in a sample-and-hold fashion, see the illustration in \rfig{controlinput}. 
As shown in \rfig{controlinput}, $t_{k+1} = t_k + \sum^N _{i=1} \delta_i$ is the next transmission time when the plant sends $x(t_{k+1})$ as the new current state information, which is obtained by the self-triggered strategy proposed in the next section. Furthermore, by making use of the flexibility of selecting control samples when multiple control inputs are allowed to be transmitted (namely when $N>1$), we will provide an efficient way of how to pick up control samples to be transmitted, such that the reduction of the communication load is achieved. 
\begin{myrem}\label{rem_infinite}
Since $\kappa(x)$ is a continuous control law, applying $\kappa(x)$ over the network as a dual mode strategy would in fact require an infinite transmission bandwidth. 
One way to avoid this issue is to apply $\kappa(x)$ under sample-and-hold fashion; 
\begin{equation*}
u(t) = \kappa (x(t_k)),\ \ t\in [t_k , t_{k}+\delta_{l}], 
\end{equation*}
where the sampling time $\delta_{l}$ is constant and needs to be small enough such that asymptotic stability is still guaranteed in $x\in \Omega (\varepsilon)$; 
see \cite{ref43} for the related analysis. 

Another way would be to apply $\kappa (x)$ directly at the plant as a stand-alone to stabilize the system, without needing any communication with the controller as soon as $x$ enters $\Omega (\varepsilon)$. This situation could be the case when the computation of $\kappa(x)$ is possible locally at the plant, while at the same time it is only feasible to solve the OCPs through the networked controller due to computational limitations.
For this case, $\kappa (x)$ does not need to be discretized since no communication is required locally at the plant. 
\end{myrem}

\section{Deriving self-triggered condition}
In this section we propose a self-triggered strategy for networked control systems as in \rfig{network}, under MPC with sample-and-hold controllers. 
Our self-triggered strategy will be derived based on the stability analysis by taking the optimal cost $J^* (x(t_k))$ as a Lyapunov candidate. Suppose again that at $t_k$ when the OCP is solved the optimal control input and the state trajectory are given by \req{optux} and the optimal cost is $J^* (x(t_k))$. 

Denoting $\Delta_n = \sum^n _{i=1} \delta_i < T_p $ for $1\leq n \leq N$, let $x(t_k + \Delta_n)$ be the actual state when sample-and-hold controllers $\{ u^{*}(t_k), \cdots, u^* (t_k + \Delta_n ) \} $ are applied with sampling intervals $\delta_1, \cdots, \delta_n$. Moreover, let $J^* (x(t_k+\Delta_N))$ be the optimal cost obtained by solving Problem 1 based on the new current state $x(t_k +\Delta_N)$. Then, the self-triggered condition, which determines the next transmission time $t_{k+1}$, is obtained by checking if the optimal cost regarded as a Lyapunov candidate is guaranteed to decrease, i.e.,
\begin{equation*}
J^* (x(t_k +\Delta_N)) -J^* (x(t_k)) <0.
\end{equation*}
For deriving this condition more in detail, we first recap from \textit{Lemma 3} in \cite{ref36} that for a quadratic stage cost (or \textit{Theorem 2.1} in \cite{ref41} for the non-quadratic case), the following result holds:
\begin{equation}\label{earlierst}
\begin{array}{lll}
J^* (x^* (t_k + \Delta_N )) - J^* (x(t_k))\\
\leq - {\displaystyle \int}^{t_k +\Delta_N} _{t_k}\!\!\!\! F(x^* (s), u^* (s)) {\rm d}s
\end{array}
\end{equation}
where $J^* (x^* (t_k+\Delta_N))$ is the optimal cost obtained by solving Problem 1 if the current state at $t_k +\Delta_N$ is $x^* (t_k +\Delta_N)$. This means that the optimal cost would be guaranteed to decrease if the actual state \textit{followed} the optimal state trajectory $x(s) = x^*(s)$ for $s\in [t_k, t_k +\Delta_N]$. 
From \req{earlierst}, we obtain
\begin{equation}\label{upbound}
\begin{array}{lll}
J^* ( {x} (t_k + \Delta_N) ) - J^* (x(t_k) ) \\
\leq  J^* ({x} (t_k + \Delta_N)) - J^* (x^{*} (t_k +\Delta_N) )\\
\ \ \  - {\displaystyle \int}^{t_k + \Delta_N} _{t_k} F({x}^* (s), u^* (s)) {\rm d}s
\end{array}
\end{equation}
where $F({x}^* (s), u^* (s))$ is known at $t_k$ when the OCP is solved. 
\begin{myrem}[Feasibility of Problem 1]\label{rem_local}
In order to obtain the stability property given by \req{earlierst}, one can see that the feasibility of Problem~ 1 needs to be guaranteed, see e.g., \cite{ref36}. Regarding establishing the feasibility of Problem~1, the existing procedures of event-triggered MPC (see e.g., \cite{ref23_2}) or periodic MPC (see e.g., \cite{ref36}) can be utilized; we can consider a feasible controller candidate given by $\bar{u}(s) = u^* (s)$ for all $s\in [t_{k+1}, t_k +T_p ]$ and $\kappa (\bar{x}(s))$ for all $s \in ( t_k +T_p, t_{k+1} + T_p ]$, to obtain \req{earlierst}. However, compared with the existing procedures, the condition $\dot{\kappa} (x) \leq K_u$ is additionally required for the existence of the local controller, such that this controller candidate becomes {admissible}. More specifically,  since we have $\dot{\kappa}(x)= \frac{\partial \kappa(x)}{\partial x} \phi(x, \kappa(x))$, $K_u$ must satisfy 
\begin{equation*}
K_u \geq  \underset{x \in \Omega (\varepsilon_f ) }{{\rm max} }\left \{ \left| \left|\frac{\partial \kappa(x)} {\partial x} \cdot  \phi(x, \kappa(x)) \right| \right|\right \},
\end{equation*}
and this needs to be computed off-line.
\end{myrem}
 For notational simplicity in the sequel, let $E_x (\delta_1, \cdots, \delta_n )$ be the upper bound of $||{x}
^* (t_k + \Delta_n ) - {x} (t_k + \Delta_n )||$ for $1\leq n \leq N$.
The following lemmas are useful to derive a more detailed expression of \req{upbound}:
\begin{mylem}\label{lem1}
Under the Assumptions $1 - 3$, the optimal cost $J^* (x)$ is Lipschitz continuous in $x\in \Sigma_V$, with Lipschitz constant $L_J$ given by
\begin{equation}\label{lipj}
L_J = \left (\cfrac{L_F}{L_{\phi}} + L_{V_f} \right ) e^{L_{\phi} T_p} - \cfrac{L_F}{L_{\phi}}
\end{equation}
\end{mylem}

For the proof of \rlem{lem1}, see Appendix. 
\begin{mylem}\label{lem2}
Suppose that the sample-and-hold controllers given by \req{sampledcontroller} 
are applied to the plant \req{sys1} from $t_k$. Then, the upper bound of $||{x}
^* (t_k + \Delta_N ) - {x} (t_k + \Delta_N )||$, denoted as $E_x (\delta_1, \cdots, \delta_{N} )$, is obtained by the following recursion for $2\leq n\leq N$:
\begin{equation}\label{erx2}
\begin{array}{lll}
E_x(\delta_1,\cdots, \delta_{n})
= E_x (\delta_1 \cdots, \delta_{n-1}) e^{L_{\phi} \delta_{n}} + h_{x} (\delta_{n})
\end{array}
\end{equation}
with $E_x (\delta_1) = h_x (\delta_1)$, where
\begin{equation}\label{hxt}
h_{x} (t) = \frac{2K_u L_G }{ {L^2 _{\phi}}}(e^{L_{\phi} t } -1 ) - \frac{2K_u L_G }{ {L_{\phi}}}t
\end{equation}
\end{mylem}
\begin{IEEEproof}
We first show $E_x (\delta_1 ) = h_x (\delta_1)$. Observe that ${x} (t_k + \delta_1 )$ and ${x}^* (t_k + \delta_1 )$ are given by
\begin{equation*}
{x} (t_k + \delta_1 ) = x(t_k ) + {\int}^{t_k + \delta_1} _{t_k} \phi ({x} (s), u^{*} (t_k)) {\rm d}s,
\end{equation*}
\begin{equation*}
{x}^* (t_k + \delta_1 ) = x(t_k ) + {\int}^{t_k + \delta_1} _{t_k} \phi ({x}^*(s), u^{*} (s)) {\rm d}s
\end{equation*}
We obtain
\begin{equation}\label{erxd}
\begin{array}{lll}
||{x}(t_k +\delta_1) -{x}^* (t_k +\delta_1) || \\
 \leq  {\displaystyle \int}^{t_k + \delta_1} _{t_k} L_{\phi} || {x}(s) - {x}^* (s) || {\rm d}s + \cfrac{1}{2} L_G K_u \delta_1 ^2
\end{array}
\end{equation}
where we have used 
\begin{equation}\label{uupperbound}
||g({x} (s)) (u^*(t_k) -u^* (s))|| \leq L_G K_u ( s-t_k)
\end{equation}
from \ras{as1} and the control input constraint $||\dot{u}(s)||\leq K_u$. Therefore, by applying the Gronwall-Bellman inequality, we obtain 
\begin{equation*}
\begin{array}{lll}
||{x}(t_k +\delta_1) -{x}^* (t_k +\delta_1) ||  \\
\leq \cfrac{2K_u L_G }{ {L^2 _{\phi}}}(e^{L_{\phi} \delta_1 } -1 ) - \cfrac{2K_u L_G }{ {L_{\phi}}}\delta_1
\end{array}
\end{equation*}
and thus $E_x (\delta_1 ) = h_x (\delta_1)$.
Now assume that $E_x (\delta_1 \cdots, \delta_{n-1})$ is given for $n\geq 2$.
We similarly obtain
\begin{equation}\label{erxd2}
\begin{array}{lll}
||{x} (t_k + \Delta_n ) - {x} ^* (t_k + \Delta_n ) || \\
\leq   || {x} (t_k + \Delta_{n-1} )- {x}^* (t_k + \Delta_{n-1} )||\\
\ \ \ + {\displaystyle \int}^{t_k +\Delta_n} _{t_k+\Delta_{n-1} } L_{\phi} ||{x}(s)-{x}^*(s)|| {\rm d}s+ \cfrac{1}{2} L_G K_u \delta_n ^2
\end{array}
\end{equation}
The only difference between \req{erxd} and \req{erxd2} is that the initial difference $||x(t_k + \Delta_{n-1} )-x^* (t_k + \Delta_{n-1} )||$ that is upper bounded by  $E_x (\delta_1, \cdots, \delta_{n-1})$ is included in \req{erxd2}. By applying the Gronwall-Bellman inequality again, we obtain
\begin{equation*}
\begin{array}{lll}
||{x} (t_k + \Delta_n ) - {x} ^* (t_k + \Delta_n ) || \\
\leq  E_x (\delta_1 \cdots, \delta_{n-1}) e^{L_{\phi} \delta_{n}} + \cfrac{2K_u L_G }{ {L^2 _{\phi}}}(e^{L_{\phi} \delta_n } -1 ) \\
\ \ - \cfrac{2K_u L_G }{ {L_{\phi}}}\delta_n
\end{array}
\end{equation*}
Thus \req{erx2} holds. 
Therefore, the upper bound $E_x (\delta_1, \cdots, \delta_N)$ is obtained by using $E_x(\delta_1) = h_x (\delta_1)$ at first, and then recursively using \req{erx2} for $n=2, \cdots, N$. This completes the proof. 
\end{IEEEproof}
Using \rlem{lem1} and \rlem{lem2}, \req{upbound} is rewritten by 
\begin{equation*}
\begin{array}{lll}
J^* ( {x}(t_k + \Delta_N) ) - J^* (x(t_k) ) \\
\leq  L_J E_x (\delta_1, \cdots, \delta_N) 
 - {\displaystyle \int}^{t_k + \Delta_N} _{t_k} F( {x}^{*} (s) , u^{*}(s) ) {\rm d}s.
\end{array}
\end{equation*}
Therefore, letting
\begin{equation}\label{trigcondmultiple}
\begin{array}{lll}
E_x(\delta_1, \cdots, \delta_N) < \cfrac{\sigma}{L_J} {\displaystyle \int}^{t_k +\Delta_N}  _{t_k}F({x}^{*} (s), u^{*}(s) ) {\rm d}s
\end{array}
\end{equation}
where $0 <\sigma <1$, we obtain
\begin{equation}\label{costdecreasing}
\begin{array}{lll}
J^* ( {x} (t_k + \Delta_N) ) - J^* (x(t_k) ) \\
\ \ \ \ \ \ \ \ \ \ \ <  (\sigma -1) {\displaystyle \int}^{t_k + \Delta_N} _{t_k} F( {x}^{*} (s) , u^{*}(s) ) {\rm d}s \\
\ \ \ \ \ \ \ \ \ \ \ < 0
\end{array}
\end{equation}
and the cost is guaranteed to decrease. 
In our proposed self-triggered MPC strategy, therefore, the next transmission time $t_{k+1}$ is determined by the time when the violation of \req{trigcondmultiple} takes place, i.e., 
\begin{equation}\label{nexttime}
t_{k+1} = {\rm inf} \left \{ \hat{t}_{k+1}  \ |\  \hat{t}_{k+1} > t_k, \Gamma({\delta} _1,\cdots, {\delta} _N) = 0 \right \},
\end{equation}
where $\hat{t}_{k+1} = t_k + \sum^{N} _{i=1} \delta_i$ and $\Gamma({\delta} _1,\cdots, {\delta} _N)$ is given by 
\begin{equation*}
\begin{array}{lll}
\Gamma (\delta_1, \delta_2, \cdots, \delta_N) \\
\ = E_x(\delta_1, \cdots, \delta_N) - \cfrac{\sigma}{L_J} {\displaystyle \int}^{\hat{t}_{k+1}}  _{t_k}F({x}^{*} (s), u^{*}(s) ) {\rm d}s.
\end{array}
\end{equation*}
Note that between $t_k$ and ${t}_{k+1}$, there exists an infinite number of patterns for the selection of sampling time intervals $\delta_1, \cdots, \delta_N$. Since $E_x (\delta_1,\cdots, \delta_N)$ in the left-hand-side (L.H.S) of \req{trigcondmultiple} depends on these intervals, the way to select $\delta_1,\cdots, \delta_N$ clearly affects the next transmission time $t_{k+1}$ obtained by \req{nexttime}. 
In the next section, we propose a way to adaptively select $\delta_1, \cdots, \delta_N$, such that the communication load can be reduced as much as possible. 
\section{Choosing sampling intervals}
In this section we provide an efficient way of adaptively selecting sampling intervals $\delta_1, \delta_2, \cdots, \delta_N$, aiming at reducing the communication load for networked control systems. 
In the following, we let $\delta^* _1, \delta^* _2, \cdots, \delta^* _N$ be the selected sampling intervals by the controller to transmit corresponding optimal control samples. 
In order to satisfy \req{trigcondmultiple} as long as possible, one may select the intervals $\delta^* _1,\cdots, \delta^* _N$ such that $E_x( \delta_1,\cdots, \delta_N)$ is minimized. This is formulated as follows:
\begin{equation}\label{opttk1}
t_{k+1} = {\rm inf} \{ \hat{t}_{k+1} \ |\ \hat{t}_{k+1} > t_k, \ \Gamma({\delta}^* _1,\cdots, {\delta}^* _N) = 0 \},
\end{equation}
where $\hat{t}_{k+1} = t_k +{\sum}^{N} _{i=1} {\delta}^* _i$ and ${\delta}^* _1,\cdots, {\delta}^* _N$ are optimal sampling time intervals between $t_k$ and $\hat{t}_{k+1}$ such that $E_x( \delta_1,\cdots, \delta_N)$ is minimized: 
\begin{equation}\label{simplest}
\begin{array}{lll}
{\delta}^* _n = \underset{\delta_1, \delta_2, \cdots, \delta_N }{\rm argmin}\ E_x(\delta_1, \cdots, \delta_N),\ \  n=1,\cdots, N \\
 \ \ \ \ {\rm s.t.}\  \hat{t}_{k+1} =t_k+{\displaystyle {\sum}^{N} _{i=1}} \delta _i.
\end{array}
\end{equation}

In this approach, it is required to solve the optimization problem \req{simplest} for each $\hat{t}_{k+1}$ and check if the self-triggered condition \req{trigcondmultiple} is satisfied. 
This means that the controller needs to both solve \req{simplest} and check \req{trigcondmultiple} until the violation $\Gamma({\delta}^* _1,\cdots, {\delta}^* _N) = 0$ occurs. 
Therefore, trying to obtain \req{opttk1} is in fact not practical from a computational point, since the optimization problem \req{simplest} needs to be solved for a possibly large number of times. 
Moreover, since the solution to \req{simplest} does not provide an explicit solution, numerical calculations of solving \req{simplest} would become more complex as $N$ becomes larger. 

Therefore, we propose a following alternative algorithm to make the problem of searching for the sampling intervals easier. In contrast to the above approach, this scheme requires only $N$ local optimizations to obtain the sampling intervals, and furthermore, a more explicit solution can be found. 

\begin{figure}[t]
   \centering
    \subfigure[Step 1: Assume $u^{*}(t_k)$ is applied, and find $\tau_1$ when \req{trigcondmultiple} is violated. ]
      {\includegraphics[width=3.7cm]{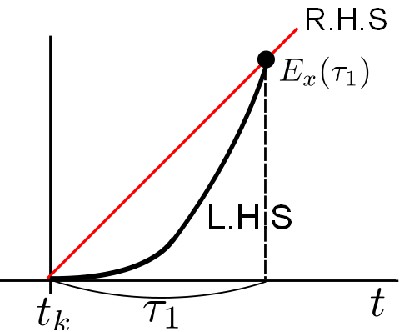}}
    \subfigure[Step 2-1: Find $0< \delta^*  _1<\tau_1$ maximizing the difference $E_x(\tau_1 ) - E_x(\delta _1, \tau_1-\delta _1)$.]
      {\includegraphics[width=4.5cm]{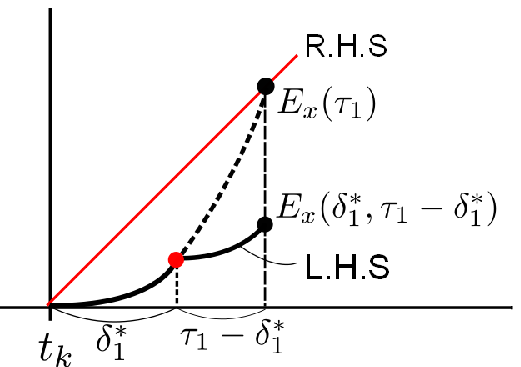}}
    \subfigure[Step 2-2: We can continue to use $u^*(t_k+\delta^* _1)$ to find the time interval $\tau_2$ until \req{trigcondmultiple} is violated.]
      {\includegraphics[width=3.6cm]{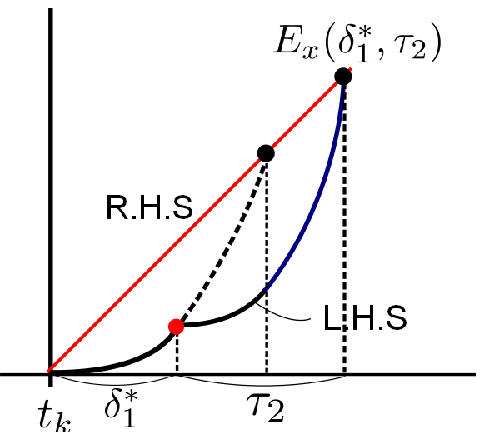}}
    \subfigure[Similarly to (b), find $0<\delta^* _2<\tau_2$ maximizing the difference $E_x(\delta^* _1,\tau_2 ) - E_x(\delta^* _1, \delta _2, \tau_2-\delta_2)$ and follow the steps until we obtain $N$ samples.]
      {\includegraphics[width=5.1cm]{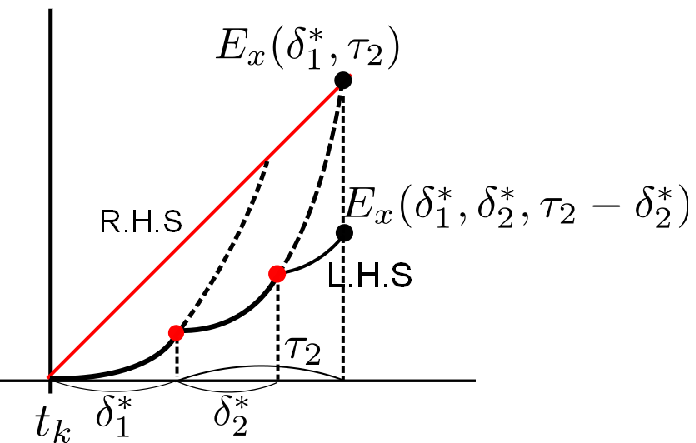}}
    \caption{The way to find sampling intervals: the L.H.S and the R.H.S are the evolutions of left-hand-side and right-hand side in \req{trigcondmultiple}.}
    \label{steps}
 \end{figure}

\noindent
\textit{Algorithm 1 (Choosing sampling time intervals)}: 
\begin{enumerate}
\item Suppose that only $u^{*}(t_k)$ is applied for $t\geq t_k$ as a constant controller, and find the time $t_k +\tau_1$ when the triggering condition \req{trigcondmultiple} is violated, see \rfig{steps} (a). We obtain $E_x(\tau_1)$ as the upper bound of $||x^* (t_k +\tau_1)-x(t_k +\tau_1)||$. 
If $N=1$, we set $\delta^* _1 = \tau_1 $.
\item If $N\geq 2$, we set $\delta^* _1\in [0, \tau_1]$ in the following way. Suppose that $u^* (t_k)$ and $u^* (t_k +\delta_1)$ are applied for [$t_k, t_k +\delta_1$], [$t_k +\delta_1, t_k +\tau_1$] respectively. 
This means we obtain $E_x (\delta_1, \tau_1 - \delta_1)$ as the upper bound of $||x(t_k + \tau_1) -x^* (t_k +\tau_1)||$.
Then, find $\delta^* _1 \in [0, \tau_1]$ which maximizes the difference of two upper bounds, i.e.,
\begin{equation*}
\delta^* _1 = \underset{\delta_1 \in [0, \tau_1]}{\rm argmax}\ \{ E_x (\tau_1 ) -E_x (\delta_1, \tau_1 -\delta_1) \},
\end{equation*}
see \rfig{steps} (b). 
As shown in \rfig{steps} (c), by maximizing the above difference, $u^*(t_k +\delta^* _1)$ can continue to be applied until the time when \req{trigcondmultiple} is again violated. We denote $\tau_2$ as the time interval when the violation of \req{trigcondmultiple} takes place after the time $t_k + \delta^* _1$. 
If $N=2$, we set $ \delta^* _2 = \tau_2$.
\item We follow the above steps until we get $N$ intervals. That is, given $n-1$ sampling intervals $\delta^* _1, \cdots, \delta^* _{n-1}$ for $2\leq n < N$, find $\tau_n$ when the triggering condition is violated to obtain $E_x(\delta^* _1, \cdots, \delta^* _{n-1}, \tau_n)$. Then, find $\delta^* _n \in [0, \tau_n]$ maximizing 
$E_x (\delta^* _1,\cdots, \delta^* _{n-1}, \tau_n) - E_x(\delta^* _1, \cdots, \delta^* _n, \tau_n - \delta^* _n)$, i.e.,
\begin{equation*}
\begin{array}{ccc}
\delta^* _n = \underset{\delta_n \in [0, \tau_n]}{\rm argmax}\ \{ E_x (\delta^* _1,\cdots, \delta^* _{n-1}, \tau_n) \\
\ \ \ \ \ \ \ \ \ \ \ \ \ \ \ \ \ \ \ \ \ \ \ \ \ \ \ \ - E_x(\delta^* _1, \cdots, \delta^* _n, \tau_n - \delta^* _n) \}.
\end{array}
\end{equation*}
For the last step at $n=N$, we set $\delta^* _N = \tau_N $, as the final time interval.
\end{enumerate}

Instead of solving the optimization problem \req{simplest} possibly for a very large number of times, Algorithm~1 requires only $N$ local optimization problems to obtain the sampling intervals $\delta^* _1, \cdots, \delta^* _{N}$. 
Algorithm~1 may not provide the largest possible next transmission time, since it does not minimize $E_x(\delta _1, \cdots, \delta _{N})$. However, as we will see through several comparisons in simulation results presented in Section~VI, Algorithm~1 is more practical than the method to obtain \req{opttk1}, as it requires much less computation time. 
Furthermore, compared with \req{simplest} that provides no explicit solutions, 
the following lemma states that the solutions to the local optimization problems can be obtained by a simple numerical procedure. 
\begin{mylem}\label{lem3}
Given $\delta^* _1, \delta^* _2, \cdots, \delta^* _{n-1}$, and $\tau_n$ for $1\leq n < N$, the transmission interval $\delta^* _n$ maximizing \\
$E_x (\delta^* _1,\cdots, \delta^* _{n-1}, \tau_n) - E_x(\delta^* _1, \cdots, \delta^* _n, \tau_n - \delta^* _n)$ in Algorithm~1, step (iii),  is obtained by the solution to
\begin{equation}\label{opttn}
e^{L_{\phi}(\tau_n - \delta^* _n)} = \cfrac{1}{(1-L_{\phi}\delta^* _n)}
\end{equation}
Furthermore, there always exists a solution of \req{opttn} satisfying $0<\delta^* _n <\tau_n$. 
\end{mylem}
\begin{IEEEproof}
From \req{erx2}, $E_x (\delta^* _1,\cdots, \delta^* _{n-1}, \tau_n)$ is given by
\begin{equation}\label{e1}
\begin{array}{lll}
E_x (\delta^* _1,\cdots, \delta^* _{n-1}, \tau_n) \\
\ \ = E_x (\delta^* _1,\cdots, \delta^* _{n-1})e^{L_{\phi} \tau_n} + h_{x} (\tau_n)
\end{array}
\end{equation}
For $E_x(\delta^* _1, \cdots,\delta^* _{n-1}, \delta_n, \tau_n - \delta_n)$, we obtain
\begin{equation*}
\begin{array}{lll}
E_x (\delta^* _1,\cdots, \delta_{n}, \tau_n-\delta_n) \\
\ \ = E_x (\delta^* _1,\cdots, \delta^* _{n-1}, \delta_n)e^{L_{\phi} (\tau_n-\delta_n)} + h_{x} (\tau_n-\delta_n)\\
\ \ = E_x (\delta^* _1,\cdots, \delta^* _{n-1})e^{L_{\phi}\tau_n} + h_{x} (\tau_n) \\
\ \ \ \ \ - \cfrac{2 K_u L_G}{L_{\phi}}  \left (e^{L_{\phi}(\tau_n - \delta_n)} -1 \right ) \delta_n
\end{array}
\end{equation*}
Thus, we obtain
\begin{equation}\label{dif_Ex}
\begin{array}{lll}
E_x (\delta^* _1,\cdots, \delta^* _{n-1}, \tau_n) - E_x(\delta^* _1, \cdots, \tau_n - \delta_n) \\
\ \ = \cfrac{2 K_u L_G}{L_{\phi}} \delta_n \left (e^{L_{\phi}(\tau_n - \delta_n)} -1 \right ) >0
\end{array}
\end{equation}
Therefore, by differentiating \req{dif_Ex} with respect to $\delta_n$ and solving for $0$, we obtain \req{opttn}.

Now it is shown that we can always find $0<\delta^* _n <\tau_n$ satisfying \req{opttn}. 
As $\delta_n \rightarrow 0$, we get 
\begin{equation*}
e^{L_{\phi}(\tau_n-\delta_n)} >\cfrac{1}{1-L_{\phi}\delta_n}
\end{equation*}
Moreover, we obtain 
\begin{equation*}
e^{L_{\phi}(\tau_n-\delta_n)}<\cfrac{1}{1-L_{\phi}\delta_n}
\end{equation*}
as $\delta_n \rightarrow \tau_n$ if $\tau_n < 1/L_{\phi}$, or $\delta_n \rightarrow 1/L_{\phi}$ if $\tau_n > 1/L_{\phi}$. Therefore, there always exists $\delta^* _n$ satisfying $0<\delta^* _n <\tau_n$. This completes the proof.  
\end{IEEEproof}
\rlem{lem3} states that $\delta^* _n$ can be found by solving \req{opttn}, once $\tau_n$ is obtained. Note that the difference \req{dif_Ex} is positive for any $0< \delta_n <\tau_n$. 
This means that if we use larger $N$, then we obtain longer transmission intervals. 

To conclude, the over-all self-triggered algorithm, including the OCP and Algorithm~1, is now stated: \\

{\textit{Algorithm 2}}: (Self-triggered strategy via adaptive control samples selection) 
\begin{enumerate}
\item At an update time $t_k$, $k \in \mathbb{N}_{\geq 0}$, if $x(t_k) \in \Omega (\varepsilon)$, then switch to the local controller $\kappa (x)$ to stabilize the system. 
Otherwise, solve Problem 1 to obtain $u^* (s )$, ${x}^* (s )$ for all $[t_k, t_k +T_p]$.
\item For a given $N$, calculate $\delta^* _1, \delta^* _2, \cdots, \delta^* _{N}$ and obtain the next transmission time $t_{k+1} =t_k + {\sum}^{N} _{i=1} {\delta}^* _i$, according to Algorithm~1. 
Then the controller transmits the following control samples to the plant;
\begin{equation}\label{controlsamples}
\{ u^{*}(t_k), u^* (t_k +\delta^* _1), \cdots, u^* (t_k + \sum^N _{i=1} \delta^* _i ) \}.
\end{equation}
\item The plant applies \req{controlsamples} in a sample-and-hold fashion, and transmits $x(t_{k+1})$ to the controller as the new current state to solve the next OCP.  
\item $k\leftarrow k+1$ and go back to Step (i).
\end{enumerate}
\begin{myrem}[\textit{Effect of time delays}]\label{delays}
So far we have ignored time delays arising in transmissions or calculations solving OCPs. In practical applications, however, it may be important to take delays into account. A method for dealing with the delays for MPC has been proposed in the recent paper \cite{ref22}, where the authors proposed delay compensation schemes by using forward prediction, i.e., even though the delays occur, the actual state is still able to be obtained from the system model \req{sys1} (see {Eq. (11)} in \cite{ref22}). 
Note, however, that in order to compensate time delays and guarantee stability, the network delays need to be upper bounded. 
More specifically, denoting $\bar{\tau}_d$ as the total maximum time delay which could arise,  then $\bar{\tau}_{d}$ needs to satisfy $\bar{\tau}_{d} < T_p - \Delta_{N}$ so that the inter-sampling time and the delay cannot exceed the prediction horizon $T_p$. Thus, assuming that $\bar{\tau}_{d}$ is known, the condition
\begin{equation}\label{delay_compensate}
\Delta_N < T_p - \bar{\tau}_{d}, 
\end{equation}
is required in the self-triggered strategy in addition to \req{trigcondmultiple}. 
\end{myrem}
\begin{myrem}[\textit{Effect of model uncertainties}]
For simplicity reasons, we have not considered the effect of model uncertainties or disturbances. However, with a slight modification of the self-triggered condition, these effects can be taken into account. Suppose that the actual state is $x^a (t)$ and the dynamics are given by $\dot{x}^a = \phi (x^a, u) + w$ where $w$ represents the disturbance or modeling error satisfying $||w|| \leq w_{\rm max}$. 
In this case, the new upper bound of $||x^* (t_k + \Delta_N) -x^a (t_k + \Delta_N)||$, denoted as $\hat{E} _x (\delta_1, \cdots, \delta_N)$ is given by 
\begin{equation*}
\hat{E} _x (\delta_1, \cdots, \delta_N) = E_x  (\delta_1, \cdots, \delta_N) + \cfrac{w_{\rm max}}{L_{\phi}} \left( e^{L_{\phi}\Delta_N} -1 \right),
\end{equation*}
where we use Gronwall-Bellman inequality, see \cite{ref20} for the related analysis. The corresponding self-triggered condition is thus given by replacing $E_x$ with $\hat{E} _x$ in \req{trigcondmultiple}. 
Similarly to Algorithm~1, it is required to obtain $\delta^* _n$ by maximizing the difference of two upper bounds $\hat{E} _x$. 
However, we can easily see that 
\begin{equation*}
\begin{array}{lll}
\hat{E}_x (\delta^* _1,\cdots, \delta^* _{n-1}, \tau_n) - \hat{E}_x(\delta^* _1, \cdots, \tau_n - \delta_n) \\ 
\ \ \ = {E}_x (\delta^* _1,\cdots, \delta^* _{n-1}, \tau_n) - {E}_x(\delta^* _1, \cdots, \tau_n - \delta_n)
\end{array}
\end{equation*}
as the effect of the disturbance can be canceled by taking the difference of the two $\hat{E} _x$.
Thus, Algorithm~1 does not need to be modified as the way to obtain sampling time intervals is not affected. 
\end{myrem}
\begin{myrem}[\textit{On the selection of the number of control samples}]\label{selection_of_N}
From \req{dif_Ex} the difference of two upper bounds is always positive, so that more time is allowed for the self-triggered condition to be satisfied by setting a new sampling time (see the illustration in \rfig{steps}(c)). 
Thus we obtain longer transmission time intervals as $N$ is chosen larger. However, $N$ needs to be carefully chosen such that the network bandwidth limitation can be taken into account; large values of $N$ may not be allowed for the network due to narrow bandwidth. 
Moreover, even though Algorithm~2 makes efficient calculations of $N$ sampling intervals, a larger selection of $N$ means more iterations of \req{opttn}, which may induce larger network delays. 
As we have already mentioned in \rrem{delays}, the delays can be compensated. However, the allowable delays must be limited as shown in \req{delay_compensate}. 
Thus, when implementing Algorithm~2, $N$ needs to be appropriately selected such that it satisfies not only the constraint for network bandwidth but also for network delays fulfilling \req{delay_compensate}. 
\end{myrem}
\section{Stability analysis}
In this section, we establish stability under our proposed self-triggered strategy. As the first step, it is shown that if the current state $x(t_k)$ is outside of $\Omega (\varepsilon)$, there always exists a positive minimum inter-execution time for the self-triggered condition \req{trigcondmultiple}, i.e., 
there exists $\delta_{\rm min} >0 $ satisfying \req{trigcondmultiple} for all $[t_k , t_k + \delta_{\rm min}]$.
We will show this only for the case where one control sample is transmitted, 
i.e, $N=1$, since larger $N$ allows for longer transmission intervals according to \rlem{lem3} and \rrem{selection_of_N}. 

The self-triggered condition for the case $N=1$ is given by $ E_x(\delta_1)  < \frac{\sigma}{L_J} {\int^{t_k +\delta_1}  _{t_k}} F(x^{*} (s), u^{*}(s) ) {\rm d}s$, where $x^*(t_k) = x(t_k)$ and $E_x(\delta_1) = h_x (\delta_1)$. By using $F(x,u) \geq \alpha_1(||x||)$ from \ras{as3}, the condition can be replaced by
\begin{equation}\label{steq1}
\begin{array}{lll}
{\displaystyle \int^{\delta_1} _{0}}\left \{ \cfrac{\sigma}{L_J} \alpha_1(||{x}^*(t_k +\eta)||) - \cfrac{2K_u L_G}{L_{\phi}} (e^{L_{\phi}\eta } -1 ) \right \} {\rm d}\eta > 0
\end{array}
\end{equation}
where $h_x(\delta_1 )$ is included in the integral. A sufficient condition to satisfy \req{steq1} is that the integrand is positive for all $0\leq \eta \leq \delta_1$, i.e., 
\begin{equation}\label{minimum}
 \alpha_1 (||{x}^* (t_k +\eta ) ||) > \cfrac{2K_u L_G L_J}{L_{\phi}\sigma} (e^{L_{\phi} \eta } -1 )
\end{equation}
for all $0\leq \eta \leq \delta_1$. 
We will thus show that if $x(t_k) \in \Sigma_V \backslash \Omega (\varepsilon)$ there exists a positive time interval $\delta_{\rm min}>0$ satisfying \req{minimum} for all $0\leq \eta \leq \delta_{\rm min}$.

\begin{figure}[t]
  \begin{center}
   \includegraphics[width=7cm]{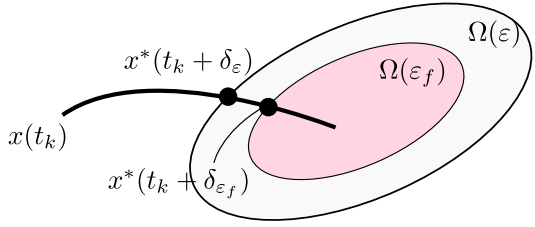}
   \caption{The illustration of $\Omega (\varepsilon)$ and the restricted terminal region $\Omega (\varepsilon_f)$.}
   \label{terregion}
  \end{center}
 \end{figure}

Suppose at a certain time $t_k +\delta_{\varepsilon}$, the optimal state ${x}^*(t_k +\delta_{\varepsilon})$ enters $\Omega (\varepsilon)$ from $x(t_k)\in \Sigma_V \backslash \Omega (\varepsilon)$, i.e., $x^* (t_k+\delta_{\varepsilon}) \in  \partial \Omega (\varepsilon )$, and it enters $\Omega (\varepsilon_f )$ at $t_k +\delta_{\varepsilon_f}$, i.e., $x^* (t_k+\delta_{\varepsilon_f}) \in  \partial \Omega (\varepsilon_f )$, as shown in \rfig{terregion}. Since $\Omega (\varepsilon_f ) \subset \Omega (\varepsilon)$, it holds that  $\delta_{\varepsilon_f}-\delta_{\varepsilon} >0$. 

To guarantee the existence of $\delta_{\rm min}$, the following two cases are considered:
\begin{enumerate}
\item ${x}^*(t_k + \eta )$ is outside of $\Omega (\varepsilon_f )$ for all the time until \req{minimum} is violated. That is,
${x}^*(t_k + \eta ) \notin \Omega (\varepsilon_f )$ for all $\eta \in [0, \bar{\eta} ]$, where 
\begin{equation}
\alpha_1 (||{x}^* (t_k +\bar{\eta} ) ||) = \cfrac{2K_u L_G L_J}{L_{\phi}\sigma} (e^{L_{\phi} \bar{\eta} } -1 )
\end{equation}
\item $x^*(t_k + \eta)$ enters $\Omega (\varepsilon_f )$ by the time \req{minimum} is violated. 
That is, there exists $\eta' \in  [0, \bar{\eta} ]$ where we obtain $x^*(t_k + \eta' ) \in \partial \Omega (\varepsilon_f )$.
\end{enumerate}

Denote $\delta_{ {\rm min}, 1}$, $\delta_{ {\rm min}, 2}$ as minimum inter-execution times for the above cases (i),  (ii), respectively. 
For the case (i), it holds that $\alpha_1 (||{x}^*(t_k+\eta )||) \geq \alpha_1(\alpha^{-1} _2( \varepsilon_f )) >0$, since we have $F(x,u) \geq \alpha_1 (||x||)$ and $V_f (x)\leq \alpha_2 (||x||) $ from \ras{as3}. 
Thus the minimum inter-execution time $\delta_{ {\rm min}, 1}$ is given by the time interval when the R.H.S in \req{minimum} reaches $\alpha_1(\alpha^{-1} _2( \varepsilon_f ))$, i.e., 
\begin{equation}\label{minimumt2} 
\delta_{{\rm min},1}  = \cfrac{1}{L_{\phi}} \ln \left (1+ \cfrac{\sigma L_{\phi} \alpha_1(\alpha^{-1} _2( \varepsilon_f )) }{2K_u L_G L_J} \right ) >0
\end{equation}
For the case of (ii), the minimum inter-execution time is $\delta_{ {\rm min}, 2}= \delta_{\varepsilon_f}-\delta_{\varepsilon}$, 
since $x(t_k)\in \Sigma_V \backslash \Omega (\varepsilon)$ and it takes at least $\delta_{\varepsilon_f}-\delta_{\varepsilon}$ for the state to reach $\Omega (\varepsilon_f )$. Thus, considering both cases, the over-all minimum inter-execution time $\delta_{\rm min}$ is positive and given by $\delta_{\rm min} = {\rm min}\ \{ \delta_{ {\rm min},1}, \delta_{ {\rm min},2} \}$. 

Based on this result, we finally obtain the following stability theorem.
\begin{mythm}
Consider the networked control system in \rfig{network} where the plant follows the dynamics given by \req{sys1}, and the proposed self-triggered strategy (Algorithm~2) is implemented. Then, if the initial state starts from $x(t_0) \in \Sigma_V \backslash \Omega (\varepsilon)$, then the state is guaranteed to enter $\Omega (\varepsilon )$ in finite time.
\end{mythm}
\begin{IEEEproof}
We prove the statement by contradiction. Starting from $x(t_0)\in \Sigma_V\backslash \Omega (\varepsilon)$, assume that the state is outside of $\Omega (\varepsilon)$ for all the time, i.e., $x(t) \in \Sigma_V\backslash \Omega (\varepsilon )$, for all $t\in [t_0, \infty )$. 

Since there exists $\delta_{\rm min}>0$, we obtain 
\begin{equation}
\begin{array}{lll}
J^*(x(t_k)) - J^*(x(t_{k-1})) \\
\ \ <   (\sigma -1){\displaystyle \int}^{t_{k}} _{t_{k-1}} F(x^* (s) , u^* (s) ){\rm d} s \\
\ \  < (\sigma -1){\displaystyle \int}^{t_{k-1}+\delta_{{\rm min}}} _{t_{k-1}} \ \ \alpha_1(\alpha^{-1} _2( \varepsilon_f )) {\rm d}s \\
\ \ = - (1- \sigma) \alpha_1(\alpha^{-1} _2( \varepsilon_f )) \ \delta_{{\rm min}} \\
\ \ = - \bar{\delta}_J < 0 
\end{array}
\end{equation}
where we denote $\bar{\delta}_J= (1- \sigma) \alpha_1(\alpha^{-1} _2( \varepsilon_f )) \ \delta_{{\rm min}}$.
Thus, we obtain
\begin{equation}\label{costv}
\begin{aligned}
J^* (x(t_k)) &-  J^*(x(t_{k-1})) < -\bar{\delta}_J \\
J^*(x(t_{k-1})) &- J^*(x(t_{k-2})) <  -\bar{\delta}_J \\
J^*(x(t_{k-2}))&- J^*(x(t_{k-3})) < -\bar{\delta}_J \\
               &\vdots   \\
J^*(x(t_1)) &- J^*(x(t_0)) < -\bar{\delta}_J
\end{aligned}
\end{equation}
Summing over both sides of \req{costv} yields 
\begin{equation}\label{vxtk}
J^*(x(t_k)) < -k \bar{\delta}_J + J^*(x(t_0)) < -k \bar{\delta}_J + J_0,
\end{equation}
where $J_0$ is defined in \rdef{def2}. This implies $J^*(t_k) \rightarrow -\infty$ as $k\rightarrow \infty$, which contradicts the fact that $J^* (x(t_k)) \geq 0$. Therefore, there exists a finite time when the state enters $\Omega (\varepsilon)$.
\end{IEEEproof}
Note again that as soon as the state reaches $\Omega (\varepsilon )$, the local control law $\kappa (x)$ is applied as a dual mode strategy. Therefore, our control objective to asymptotically stabilize the system to the origin is achieved, i.e., $x(t) \rightarrow 0$ as $t\rightarrow \infty$.

\section{Simulation examples}
In this section we illustrate our proposed self-triggered scheme for both linear and nonlinear systems. 
Simulations were implemented in MATLAB on a PC having 2.50 GHz Intel (R) Core (TM) CPU and 4.00 GB RAM. 
As a software package, we used SNOPT in order to solve (non)linear optimal control problems, see \cite{ref42}. 
\subsection{Linear case}
An interesting example is to check if we can guarantee to stabilize \textit{un}-stable systems under our proposed aperiodic control execution. 
Therefore, as one of such examples we consider the following linearized system of inverted pendulum on a cart problem (see \cite{ref5});
\begin{equation*}
\dot{x} = Ax + Bu
\end{equation*}
where we denote $x = [x_1\ x_2\ x_3\ x_4]^\mathsf{T} \in \mathbb{R}^4$,\ \ $u\in \mathbb{R}$ and
\begin{equation*}
A=\left [
\begin{array}{cccc}
0 &  1 & 0 & 0 \\
0 & 0 & -mg/M & 0 \\
0 & 0 & 0     & 1 \\
0 & 0 & g/l  & 0
\end{array}
\right ],\ \ B= \left [
\begin{array}{c}
0 \\
1/M \\
0 \\
-1/Ml
\end{array}
\right ]
\end{equation*} 
We set $m=0.55$ as the point mass, $M=15$ as the mass of the cart, and $l=9$ as the length of the massless rod. 
The system is unstable having a positive eigenvalue $1$ in matrix $A$. The constraint for the control input is assumed to be given by $||u||\leq 8.5$. The computed Lipschitz constants $L_f$ and $L_G$ are given by $L_f =1.05$, $L_G = 0.067$. The stage and the terminal cost are assumed to be quadratic and given by $F(x, u ) = x^\mathsf{T} Q x + u^\mathsf{T} R u$ where $Q= 3.0 I_4$ and $R= 1.5$. 
\begin{figure}[t]
   \centering
   \subfigure[State trajectories of $x_1$ and $x_2$ ]
      {\includegraphics[width=8.0cm]{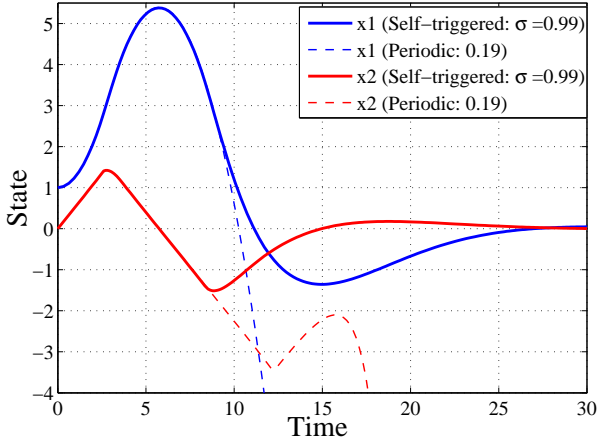}}
         \hspace{2pt}
    \subfigure[State trajectories of $x_3$ and $x_4$]
      {\includegraphics[width=8.0cm]{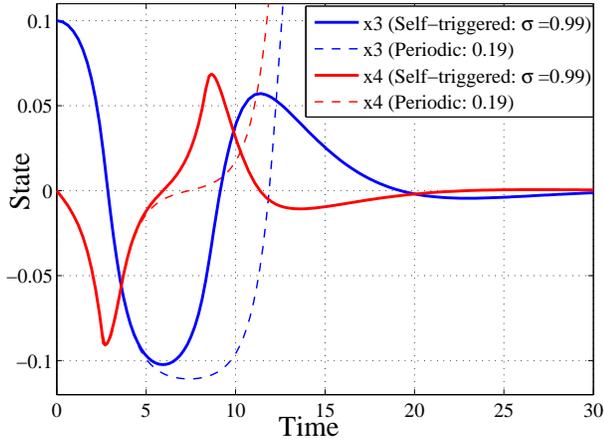}}
    \caption{State trajectories under the self-triggered MPC with $\sigma =0.99$ (solid lines) and  periodic MPC with the same average transmission interval $0.19$ (dot lines). }
    \label{sim_lin_x12}
 \end{figure}
\begin{figure}[htbp]
   \centering
   \subfigure[State trajectories of $x_1$ and $x_2$]
      {\includegraphics[width=8.0cm]{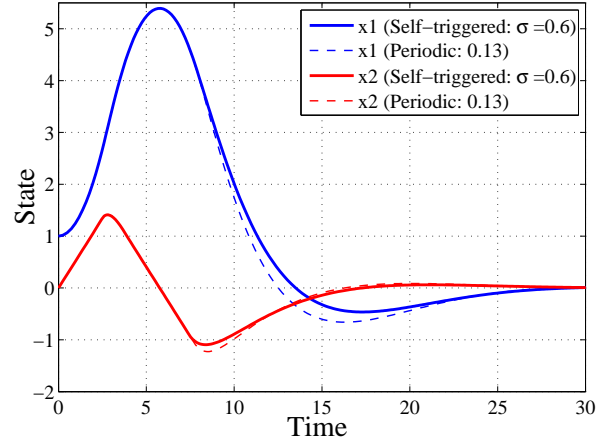}}
         \hspace{2pt}
    \subfigure[State trajectories of $x_3$ and $x_4$]
      {\includegraphics[width=8.0cm]{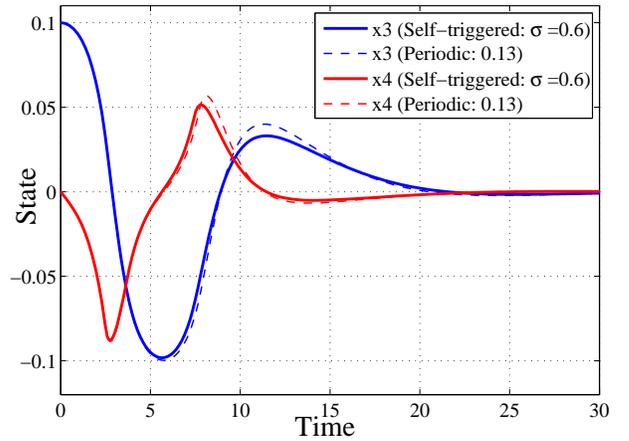}}
    \caption{State trajectories under the self-triggered MPC with $\sigma =0.6$ (solid lines) and  periodic MPC with the same average transmission interval $0.13$ (dot lines).}
    \label{sim_lin_x1206}
 \end{figure}
The terminal cost is given by $V_f = x^\mathsf{T} P_f x$, where 
\begin{equation}
P_f =\left [
\begin{array}{cccc}
21 & 70  & 490 & 507 \\
70 & 374 & 2866 & 2968 \\
490 & 2866 & 39212     & 40002 \\
507 & 2968 & 40002  & 40822
\end{array}
\right ].
\end{equation}
The matrix $P_f$ and the local terminal controller $\kappa(x) = K x$ are obtained by the procedure presented in \cite{ref36}, and given by 
\begin{equation*}
K = \left [
\begin{array}{cccc}
-1.414  & -9.739 &  -228.1& -230.9
\end{array}
\right ]
\end{equation*} 
and $\varepsilon=0.43$. We set $\varepsilon_f = 0.2$ and $K_u = 2.0$. The prediction horizon is $T_p = 14$ and the number of control input is simply given by $N=1$. 
The initial state is assumed to be $x_0 = [1\ 0\ 0.1\ 0]^\mathsf{T}$. From \req{costdecreasing}, $\sigma$ is the parameter that restricts how much the optimal cost $J^*$ is guaranteed to decrease. 
Thus we consider two cases for the choice of $\sigma$; 
$\sigma =0.6$ and $0.99$ in order to compare the control performance. 

\begin{figure}[t]
  \begin{center}
   \includegraphics[width=8cm]{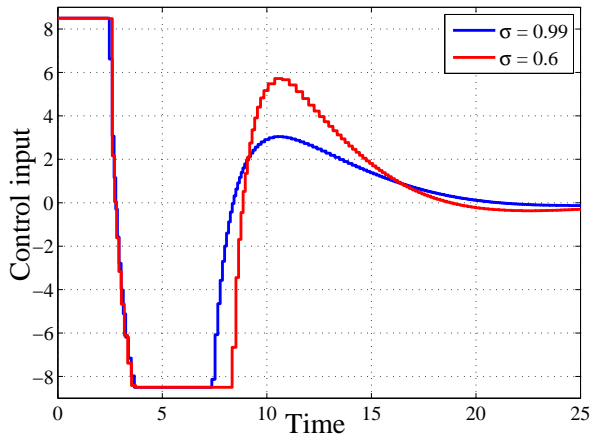}
   \caption{Control input $u$}
   \label{sim_lin_controlinput}
  \end{center}
 \end{figure}

\begin{figure}[t]
  \begin{center}
   \includegraphics[width=8cm]{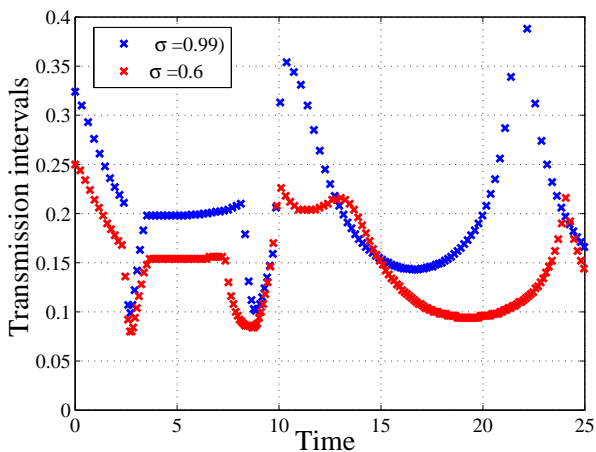}
   \caption{Transmission interval}
   \label{sim_lin_trig}
  \end{center}
 \end{figure}

\begin{figure}[t]
  \begin{center}
   \includegraphics[width=8cm]{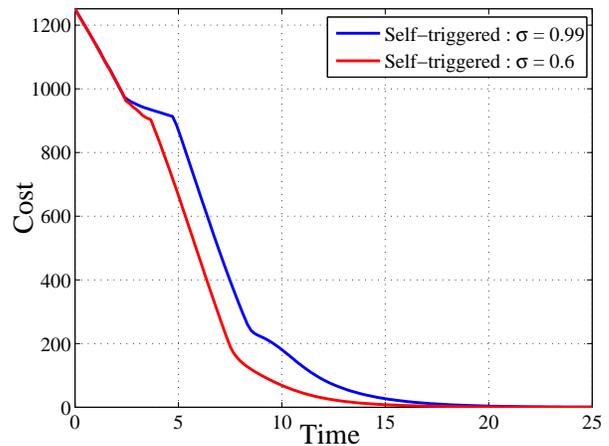}
   \caption{Cost sequence}
   \label{sim_lin_cost}
  \end{center}
 \end{figure}

\rfig{sim_lin_x12} and \rfig{sim_lin_x1206} show the state trajectories under Algorithm~2 with $\sigma=0.6$ and $\sigma=0.99$. \rfig{sim_lin_controlinput} shows the control input $u$. \rfig{sim_lin_trig} plots the transmission intervals at each update time until the state reaches $\Omega (\varepsilon )$. \rtab{av} shows the average transmission intervals for both $\sigma=0.6$ and $\sigma =0.99$. 
\rfig{sim_lin_cost} shows the sequence of the optimal costs obtained by solving the OCPs. 

As shown in Fig. 6, 7 and 8, the system is stabilized to the origin while satisfying the input constraint $||u(t)|| \leq 8.5$. Furthermore, as shown in \rtab{av} the average transmission interval becomes smaller for the case $\sigma =0.6$ than for the case $\sigma =0.99$, meaning that it requires more transmissions for smaller choice of $\sigma$. 
As for the result of optimal costs, on the other hand,  selecting smaller $\sigma$ leads to better control performances since it enforces the optimal cost to decrease more, see \rfig{sim_lin_cost}. Thus the result implies the trade-off between obtaining control performance and the transmission rate; selecting larger $\sigma$ leads to smaller number of transmissions, but it degrades the control performances. 
On the other hand, selecting smaller $\sigma$ leads to better control performances but requires more transmissions. 

To make further comparisons, we also plotted the result of state trajectories under the periodic (standard) MPC with sample-and-hold implementation in Fig. 6 and 7 as red and blue dotted lines. 
In order to compare control performances under the same transmission rate, the sampling times were selected as $0.13$ and $0.19$, which are same as the average transmission intervals obtained by Algorithm~2 in \rtab{av}. 
From \rfig{sim_lin_x1206}, the state trajectories under the self-triggered strategy with $\sigma=0.6$ have similar convergences to the periodic case. 
However, as shown in \rfig{sim_lin_x12}, the state for the periodic case with $\sigma = 0.99$ fails to be stabilized. 
This means that the sampling time was not appropriately selected as it was not chosen to guarantee stability, even though the transmission rate is the same as the self-triggered case.
Therefore, we show by considering the unstable system 
that the system is stabilized under our proposed self-triggered strategy. 
\begin{table}[t]
\begin{center}
\caption {Average transmission interval}
\label {av}
\begin{tabular}{|c|c|c|} \hline 
                           &  $\sigma =0.6$    & $\sigma =0.99$  \\ \hline \hline
Triggering interval   &    $0.13$                & $0.19$     \\ \hline
\end{tabular} 
\end{center}
\end{table}

\subsection{Nonlinear case}
For the nonlinear case, we consider the position control of a non-holonomic vehicle regulation problem in two dimensions \cite{ref39}. 
The dynamics can be modeled as
\begin{eqnarray}\label{vehicle}
\cfrac{{\rm d}}{{\rm d}t}\left[
\begin{array}{c}
x \\
y \\
\theta
\end{array}
\right] =\left[
\begin{array}{cc}
\cos \theta & 0  \\
\sin \theta & 0  \\
0 & 1
\end{array}
\right]  \left[
\begin{array}{c}
v \\
\omega
\end{array}
\right].
\end{eqnarray}

\begin{figure}[t]
  \begin{center}
   \includegraphics[width=4cm]{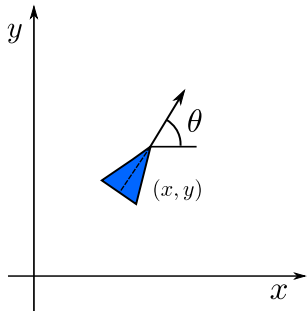}
   \caption{The state variables for a vehicle regulation problem in two dimensions.}
   \label{agent_fig}
  \end{center}
 \end{figure}

Here we denote the state as $\chi = [x\ y\ \theta]^\mathsf{T}$, consisting of the position of the vehicle $[x\ y]$, and its orientation $\theta$ (see \rfig{agent_fig}). $u=[v\ \omega]^\mathsf{T}$ is the control input and the constraints are assumed to be given by $||v||\leq \bar{v}=1.5$ and $||\omega||\leq \bar{\omega}=0.5$. 
The computed Lipschitz constant $L_{\phi}$ and a positive constant $L_G$ are given by $L_{\phi}=\sqrt{2}\bar{v}$ and $L_G =1.0$ (see \cite{ref20}). The stage and the terminal cost are given by $F= \chi^T Q \chi + u^T R u$, and $V_f = {\chi^T \chi}$ where $Q=0.1 I_3$ and $R=0.05 I_2$. The prediction horizon is $T_p =7$. Since the linearized system around the origin is uncontrollable, the procedure presented in \cite{ref38} is adopted to obtain a local controller satisfying \ras{as2}, and the parameter for characterizing the terminal set is $\varepsilon=0.8$. We set $\varepsilon_f = 0.4$ and the local controller is admissible if $K_u=1.5$.

\begin{figure}[t]
  \begin{center}
   \includegraphics[width=8cm]{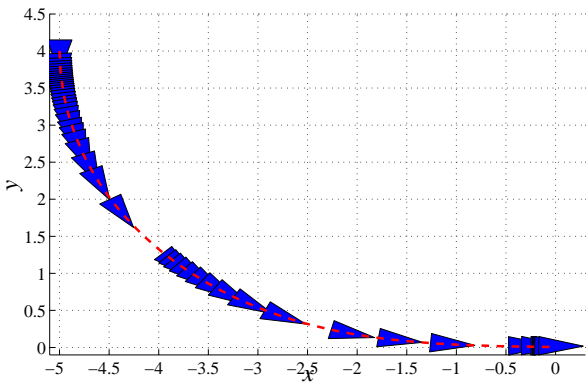}
   \caption{State trajectories under the self-triggered MPC (triangles) and periodic MPC (red dotted line). The triangles appear when the control input samples are transmitted. }
   \label{sim_non_tr}
  \end{center}
 \end{figure}

\begin{figure}[t]
  \begin{center}
   \includegraphics[width=8.0cm]{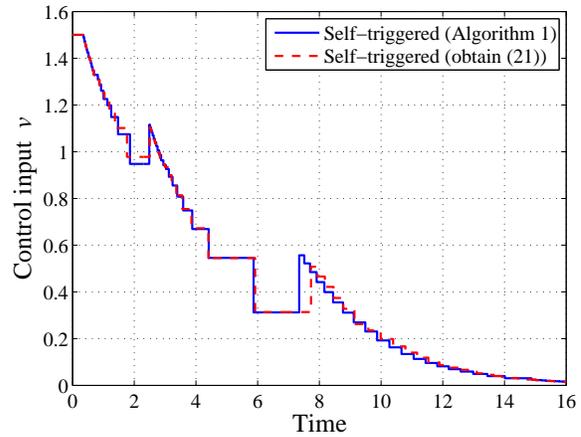}
   \caption{Control input of the self-triggered strategy by applying Algorithm~1 (blue line) and by obtaining \req{opttk1} (red dotted line). }
   \label{sim_non_input}
  \end{center}
 \end{figure}

\begin{figure}[t]
  \begin{center}
   \includegraphics[width=8.0cm]{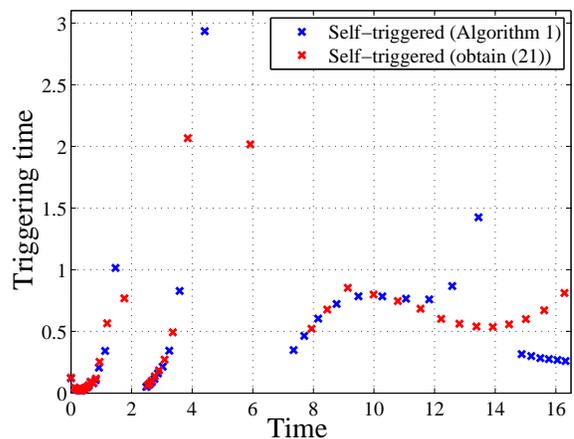}
   \caption{Transmission time intervals of the self-triggered strategy under Algorithm~1 (Blue) and the method of obtaining \req{opttk1} (Red). }
   \label{sim_non_trigtime}
  \end{center}
 \end{figure}

\rfig{sim_non_tr} shows the trajectory of the vehicle under Algorithm~2 with $\sigma =0.99$ and $N=2$, starting from the initial point $[-5\ 4\ -\pi/2]$ and its goal is the origin. The blue triangles show the position of the vehicle, where the triangle appears when control samples are transmitted to solve the OCP. The heading of the triangle shows the moving direction of the vehicle. 
\rfig{sim_non_input} shows the control input $v$. 
The triggering intervals are plotted in \rfig{sim_non_trigtime}. 
\begin{table}[t]
\centering
\caption {Average transmission intervals and the calculation times}
\subtable[Average transmission intervals]{\label{com_time_table}
\begin{tabular}{|c|c|c|c|c|c|} \hline
        $N$                    &  $1$    & $2$ & $3$ & $4$ & $5$ \\ \hline \hline
Algorithm 1              &    $0.17$  & $0.24$  & $0.41$ & $0.58$ & $0.68$ \\ \hline 
Obtain \req{opttk1}  &    $0.17$  & $0.30$  & $0.58$ & $0.74$  & $0.83$ \\ \hline 
\end{tabular}
}
\subtable[Average calculation time per each step (sec)]{\label{cal_time_table}
\begin{tabular}{|c|c|c|c|c|c|}\hline 
    $N$                     &  $1$    & $2$ & $3$ & $4$ & $5$ \\ \hline \hline
Algorithm 1              &    $0.31$  & $0.35$  & $0.41$ & $0.47$ & $0.52$ \\ \hline 
Obtain \req{opttk1}  &    $0.31$  & $7.78$  & $11.9$ & $16.3$ & $35.8$ \\ \hline 
\end{tabular}
}
\end{table}

In Section~IV, it is stated that the sampling time intervals $\delta^* _1, \cdots, \delta^* _N$ can also be obtained by solving the optimization problem \req{simplest} to obtain \req{opttk1}. To make a comparison, we also plotted the state trajectory and control input in \rfig{sim_non_tr}, \rfig{sim_non_input} as red dotted lines, where the next transmission time is given by \req{opttk1}. 
The result of transmission time intervals is also plotted as red marks in \rfig{sim_non_trigtime}. \rtab{com_time_table} shows the average transmission time intervals under both Algorithm~1 and the method of obtaining \req{opttk1} by solving \req{simplest} with different values of $N$. 
Furthermore, \rtab{cal_time_table} shows the average online computation time per each time step, including the time to solve the OCP and to select sampling time intervals $\delta^* _1, \cdots, \delta^* _N$. 
Since the minimum values of $E_x$ are obtained by solving \req{simplest}, larger average transmission time intervals are obtained than by applying Algorithm~1, as shown in \rtab{com_time_table}. 
However, as shown in \rtab{cal_time_table}, obtaining \req{opttk1} requires much more computation time, 
since \req{simplest} needs to be solved for a large number of times. 
Therefore, it is shown that Algorithm~1 is more useful for practical applications than obtaining \req{opttk1}, in terms of the calculation cost to obtain the sampling time intervals $\delta^* _1, \cdots, \delta^* _N$. 

\section{Conclusion and Future work}
We proposed an aperiodic formulation of MPC for networked control systems, where the plant with actuator and sensor systems are connected to the controller through wired or wireless sensor networks. Our proposed scheme not only provided when to solve OCPs but also the efficient way to select sampling intervals to achieve transmission intervals as large as possible. Stability under sample-and-hold implementation was also shown by guaranteeing the positive minimum inter-execution time of the self-triggered strategy. Our proposed framework was also validated through both linear and nonlinear simulation examples. Future work is to consider more detailed analysis of self-triggered strategies under additive noise or uncertainties.

\appendix
\noindent
\textit{(Proof of \rlem{lem1})}: Consider the optimal costs $J^*(x_1)$, $J^*(x_2)$ obtained by different initial states $x(0) = x_1$, $x(0) = x_2$. Here the current time is assumed to be $0$ without loss of generality.
Let $x^* _1(s), u^* _1(s)$ ($x^* _1(0)=x_1$), and $x^* _2(s), u^* _2(s)$ ($x^* _2(0)=x_2$) be the optimal state and control trajectory for $s\in [0, T_p]$, obtained by solving Problem 1. These optimal costs are then given by
\begin{equation}
J^* (x_i ) = {\displaystyle \int}^{T_p} _{0} F (x^{*} _i (s), u^{*} _i (s)) {\rm d}s + V_f (x^{*}_i (T_p))
\end{equation}
for $i=1,2$.
Now consider the difference $J^*(x_1) - J^* (x_2)$.
Assume that from the initial state ${x} _1$,
an alternative control input $\bar{u} _1(s) = u^* _2(s) \in {\cal U}$ ($s\in [0, T_p]$) is applied and let $\bar{x} _1(s)$ be the corresponding state obtained by applying $\bar{u}_1(s)$. Also let $\bar{J}(x_1)$ be the corresponding cost. Since $J^* (x_1) \leq \bar{J}(x_1)$, we obtain
\begin{equation}\label{eqV}
\begin{array}{lll}
J^*(x_1) - J^* (x_2) &\leq & {\displaystyle \int}^{T_p} _{0} L_F ||\bar{x} _1(s) - x^* _2(s)|| {\rm d}s \\
& &+ L_{V_f} ||\bar{x} (T_p) - x^* _2(T_p) ||
\end{array}
\end{equation}
where the Lipschitz continuities of $F$ and $V_f$ are used. From Gronwall-Bellman inequality, we have $||\bar{x} _1(s) - x^* _2(s)||\leq e^{L_{\phi} s} ||x_1 - x_2 ||$ for $s\in [0, T_p]$. Thus, we obtain
\begin{equation*}
\begin{array}{lll}
J^*(x_1) - J^*(x_2) \\
\ \ \ \ \leq  L_F ||x_1-x_2|| {\displaystyle \int}^{T_p} _{0} e^{L_{\phi} s} {\rm d}s + L_{V_f} e^{L_{\phi}T_p} ||x_1 - x_2 ||\\
\ \ \ \  = \left \{ \left (\cfrac{L_F}{L_{\phi}} + L_{V_f} \right ) e^{L_{\phi} T_p} - \cfrac{L_F}{L_{\phi}} \right \} ||x_1 -x_2||
\end{array}
\end{equation*}
Thus the proof is complete.

\end{document}